\documentclass[10pt, twoside]{article}
\usepackage{amsmath,amssymb,wasysym}
\usepackage{graphicx}
\usepackage[left=3cm,right=3cm, top=3.3cm, bottom=3cm]{geometry}
\newtheorem{theo}{Theorem}[section]

\newtheorem{coro}[theo]{Corolary}

\newtheorem{defi}[theo]{Definition}
\newtheorem{remark}[theo]{Remark}
\def\proof {{\noindent \bf{Proof:\hspace{4pt}}}}
\def\endproof{\hfill$\square$\vspace{6pt}}
\numberwithin{equation}{section}
\title{
{\bf\Large  Nonlinear Stability for the Periodic and Non-Periodic Zakharov System }}
\author{{\bf\large Jaime Angulo Pava}\footnote{Email: angulo@ime.usp.br}\hspace{2mm}
{\bf\large}\vspace{1mm}\\
{\it\small Department of Mathematics, IME-USP}\\
 {\it\small Rua do Mat\~ao 1010, Cidade Universit\'aria,}\\
{\it\small  CEP 05508-090, S\~ao Paulo, SP, Brazil}\vspace{3mm}\\
{\bf\large Carlos Banquet Brango}\footnote{Email: cbanquet@sinu.unicordoba.edu.co}\hspace{2mm}
{\bf\large}\vspace{1mm}\\
{\it\small Departamento de Matem\'aticas y Estad\'istica}\\
{\it\small  Universidad de C\'ordoba}\\
{\it\small Carrera 6 No. 76-103, Monter\'ia, C\'ordoba, Colombia}\vspace{3mm}}
\date{}
\begin{document}
\maketitle
\begin{abstract}
We prove the existence of a smooth curve of periodic traveling wave solutions for the Zakharov system. We also show that this type of solutions are nonlinear stable by the periodic flow generated for the system mentioned before. An improvement of the work of Ya Ping \cite{Wu1} is made, we prove the stability of the solitary wave solutions associated to the Zakharov system. \\ 


\end{abstract}

\section{Introduction}
In this essay we study the periodic Zakharov system
\begin{equation} \label{equaZakha}
\left \{
\begin{aligned}
iu_{t}+u_{xx}&=uv\\
v_{tt}-v_{xx}&=(|u|^2)_{xx},\\
\end{aligned} \right.
\end{equation}
where $u=u(x,t)\in\mathbb{C},\ v=v(x,t)\in\mathbb{R}$ and $x,t\in\mathbb{R}.$ This system was introduced by Zakharov in \cite{Zakharov1} to describe the long wave Langmuir turbulence in a plasma. The function $u=u(x,t)$ represents the slowly varying envelope of the highly oscillatory electric field and  $v$ denotes the deviation of the ion density from the equilibrium.\\

The goal of this paper is to establish the existence and nonlinear stability of periodic traveling wave solutions for the Zakharov system. More precisely, we are interested in solutions for (\ref{equaZakha}) of the form
\begin{equation}\label{solforms}
u(x,t)=e^{-i\omega t}e^{i\frac{c}{2}(x-ct)}\phi_{\omega,c}(x-ct)\ \ \ \text{and}\ \ \ 
v(x,t)=\psi_{\omega,c}(x-ct),
\end{equation}
where $\omega,c\in\mathbb{R}$ and $\phi_{\omega,c},\psi_{\omega,c}: \mathbb{R}\rightarrow\mathbb{R}$ are periodic smooth functions with the same fundamental period $L>0.$ As far as we know any result of stability for this type of waves has been established before. The first work about existence and nonlinear stability of periodic waves was made by Benjamin in \cite{benjamin3}, where he studied  periodic waves of cnoidal type for the Korteweg-de Vries equation. This work had some gaps on central parts of the stability theory that was revised and complemented by Angulo, Bona and Scialom in \cite{anguloBonaScialom}. In the last few years some papers about  the nonlinear stability on the periodic case have appeared in the literature, see for instance \cite{AnguloLibro, angulo5, angulo4, AnguloNatali2, anguloNatali, GallayHaragus1, GallayHaragus2, haragus1, NataliPastor, Neves1}. \\

Substituting the type of solutions given in (\ref{solforms}) in the system (\ref{equaZakha}), we get that $\phi=\phi_{\omega,c}$ and $\psi=\psi_{\omega,c}$ have to satisfy the next system of ordinary differential equations,
\begin{equation}\label{systemedo}
\left \{
\begin{aligned}
&(c^2 -1)\psi ''=(\phi^2)''\\
&\phi ''+\left(\omega+\frac{c^2}{4}\right)\phi=\phi\psi.
\end{aligned} \right.
\end{equation}
Integrating the first equation of the system (\ref{systemedo}) and substituting on the second one, we obtain after some algebra that the solution $\phi$ has to satisfy
\[\left(\phi'\right)^2=\frac{1}{2(1-c^2)}F(\phi),\]
where $F$ is the polynomial given by
\[F(t)=-t^4+2(1-c^2)\left(-\omega-\frac{c^2}{4}\right)t^2+4(1-c^2)A_{\phi}\]
and $A_{\phi}$ is a constant of integration. It is clear that the solutions of the equation (\ref{equaZakha}) depend of the roots of the polynomial $F.$ Assuming that $F$ has roots $\pm\eta_1$ and $\pm\eta_2$ with $0<\eta_2<\eta_1,$ we obtain the smooth curve of dnoidal waves
\[\nu\in\left(\frac{2\pi^2}{L^2},+\infty\right)\longmapsto\left(\psi_{\nu},\phi_{\nu}\right)\in H_{per}^n([0,L])\times H_{per}^n([0,L]),\ \ \text{for all}\ \ n\in \mathbb{N}, \]
with $\phi_{\nu}$ and $\psi_{\nu}$ given by
\[\phi_{\nu}(\xi)=\eta_1\text{dn}\left(\tfrac{\eta_1\xi}{\sqrt{2(1-c^2)}} ;k\right)\ \ \text{and}\ \ \ \psi_{\nu}(\xi)=-\frac{\eta^2_1}{1-c^2}\text{dn}^2 \left(\tfrac{\eta_1\xi}{\sqrt{2(1-c^2)}};k\right).\] Here,  $k^2=\frac{\eta^2_1-\eta^2_2}{\eta^2_1},$ $\nu=-\left(\omega+\frac{c^2}{4}\right)$ and dn denotes the Jacobi elliptic function of dnoidal type. This solutions are constructed  with the same fixed minimal period $L>0,$ not necessarily large.\\

With respect to the well-posedness problem for the Zakharov system, on the periodic case, this was studied by Bourgain in \cite{Bourgain2}, where a global well-posedness result was obtained for initial data $(u(0), v(0),v_t(0) )\in H^1_{per}\times L^2_{per}\times H^{-1}_{per}.$ It is worth to note that in the periodic case there exists another \textit{more general} result about well-posedness for the Zakharov system obtained by Takaoka in \cite{Takaoka1}, but for our purpose the result established by Bourgain is good enough. See also Guo and Shen \cite{GuoBoling1}, where the existence of classical periodic solutions for the system (\ref{equaZakha}) is proved. On the continuous case the Cauchy problem associated to the Zakharov system in one and several dimensions have been studied extensively, see for instance \cite{AddedAdded1, bejenaru1, BourgCollia1, Colliander1, GiniTsutVelo1, KenigPonceVega1, OzawaTsut1, Pecher, SchotetWeiste1, SulemSulem1}.\\

In order to establish the spectral properties of some linear operators which appear in the proof of the stability, we use the Floquet theory, more precisely we use the Oscillation Theorem (see Magnus and Winkler \cite{magnus}). Our spectral analysis depends basically of the next periodic and semi-periodic eigenvalue problems  associated to the Lam\'e equation, given respectively by
\[\left\{
\begin{aligned}
y''+&[\lambda-m(m+1)k^2\text{sn}^2(x,k)]y=0\\
y(0)&=y(2K(k)),\ \ y'(0)=y'(2K(k))
\end{aligned} \right.\]
and
\[\left\{
\begin{aligned}
y''+&[\lambda-m(m+1)k^2\text{sn}^2(x,k)]y=0\\
y(0)&=-y(2K(k)),\ \ y'(0)=-y'(2K(k)),
\end{aligned} \right.\]
where $\lambda\in\mathbb{R},$  $m\in\mathbb{N},$ sn denotes the Jacobi elliptic function of snoidal type and $K$ is the complete elliptic integral of the first type (see Byrd and Friedman \cite{byrdFriedman}). Recently, Neves in \cite{Neves1} proved that is possible to characterize the eigenvalues of the Hill operator $L(y)=-y''+Q(x)y$ in $L^2[0,\pi]$ if we know explicitly one of the eigenfunctions associated to this eigenvalue (in this case, $Q$ is a $C^2$ periodic function with period $\pi$). Unfortunately, we only had access to this work when we already had concluded our spectral results using the associated Lam\'e equation. We are completely sure that this new theory can be use to obtain the spectral properties of the operator studied in this paper.\\

To obtain our result of stability for the dnoidal wave solutions, we rewrite the Zakharov system as
\begin{equation}\label{ZakNovoInt}
\left \{
\begin{aligned}
v_t &=-V_x, \ \int_0^L V(x,t)dx  = 0 \\
 V_t &= -(v+|u|^2)_x \\
iu_t &+ u_{xx} =uv
\end{aligned} \right.
\end{equation}
and we adapt  to the periodic case the ideas established by Benjamin \cite{benjamin1}, Bona \cite{bona2} and Weinstein \cite{weinstein3}, then we impose the restriction
\[\int_0^L v_0(x)dx\leq \int_0^L \psi_{\nu}(x)dx,\]
where $v(x,0)=v_0(x),$ to obtain that the dnoidal waves with $c\in(-1,1)$ fixed and $\nu>\frac{2\pi^2}{L^2},$ are orbitally stable in
\[X:=L^2_{per}([0,L])\times\widetilde{L}^2_{per}([0,L])\times H^1_{per}([0,L])\]
by the periodic flow of the system (\ref{ZakNovoInt}). Here,  $\widetilde{L}^2_{per}$ is given by
\[\widetilde{L}^2_{per}([0,L])=\left\{f \in L^2_{per}([0,L]): \int_0^L f(x) dx=0\right\}.\]

With regard to the existence and stability of solitary wave solutions for the Zakharov system, there exists a result obtained by Ya Ping in \cite{Wu1}, this work is not completely right. In \cite{Wu1} the author considered the \textit{equivalent} system 
\[
\left \{
\begin{aligned}
v_t&=V_{xx},\\
V_t&=v+|u|^2,\\
iu_{t}&+u_{xx}=uv.\\
\end{aligned} \right.
\]
 Therefore, the solitary wave solution $V(x,t)=\varphi_{\omega,c}(x-ct)$ is given by
\begin{equation}\label{solWuvarphi}
\varphi_{\omega,c}(\xi)=c\sqrt{-4\omega-c^2}\tanh\left(\frac{\sqrt{-4\omega-c^2}}{2}\xi\right).
\end{equation}
Observe that this solution is not in any Sobolev space $H^s(\mathbb{R})$ and Ya ping proved stability in $L^2(\mathbb{R})\times H^1(\mathbb{R})\times H^1(\mathbb{R}),$ which is not right because the solution (\ref{solWuvarphi}) is not in the space where  the author proves the stability.  One of the goals of this paper is to improve the result of stability, for the solitary waves solutions, obtained by Ya Ping. Following the ideas  used to establish the stability on the periodic case we prove that the solitary wave solutions 
\[
\psi_{\omega,c}(\xi)=\left(2\omega+\frac{c^2}{2}\right)\ \text{sech}^2\left(\frac{\sqrt{-4\omega-c^2}}{2}\xi\right), \ \ \ \phi_{\omega,c}(\xi)=\sqrt{\frac{(-4\omega-c^2)(1-c^2)}{2}}\ \text{sech}\left(\frac{\sqrt{-4\omega-c^2}}{2}\xi\right)
\]
\[
\text{and}\ \ \ \varphi_{\omega,c}(\xi)=c\left(2\omega+\frac{c^2}{2}\right)\ \text{sech}^2\left(\frac{\sqrt{-4\omega-c^2}}{2}\xi\right)
\]
are orbitally  stable in $X=L^2(\mathbb{R})\times L^2(\mathbb{R})\times H^1(\mathbb{R})$ by the flow generated by the Zakharov system if $1-c>0$ and $4\omega+c^2\geq 0.$  It is worth to note that in the continuous case the restriction imposed above for the initial datum $v_0$  is not necessary, because using the property that the solitary wave solutions converges to zero, when $\xi$ goes to infinity, the term that force to impose this condition disappears.\\

The plan of the paper is as follows. The next section is devoted to describe briefly the notation that will be used,
and to make a few preliminary remarks regarding periodic and nonperiodic Sobolev spaces. In Section 3 we prove the existence of a smooth curve of dnoidal wave  solutions for the system (\ref{equaZakha}).  Section 4 contains the spectral analysis of some linear operators necessary to obtain our result of stability.  In Section 5 we present the result of nonlinear stability for the  dnoidal wave solutions  of the system (\ref{equaZakha}). Finally, in Section 6 we present the result of stability of the solitary waves associated to the Zakharov system.
\section{Notation}
The $L^2$-based Sobolev spaces of periodic functions are defined as follows (for further details see Iorio and Iorio \cite{ioriolibro}). Let $\mathcal{P}=C^{\infty}_{per}$ denote the collection of all functions $f:\mathbb{R}\rightarrow \mathbb{C}$ which are $C^{\infty}$ and periodic with period $2L>0.$ The collection $\mathcal{P}'$ of all continuous linear functionals from $\mathcal{P}$ into $\mathbb{C}$ is the set of \textit{periodic distributions.} If $\Psi\in \mathcal{P}'$ then we denote the value of $\Psi$ at $\varphi$
by $\Psi(\varphi)=\langle\Psi,\varphi\rangle.$ Define the functions $\Theta_k(x)=\exp(\pi ikx/L), \ k\in \mathbb{Z},\ x\in\mathbb{R}.$ The Fourier transform of $\Psi$ is the function $\widehat{\Psi}:\mathbb{Z}\rightarrow\mathbb{C}$ defined by the formula $\widehat{\Psi}(k)=\frac 1{2L}\langle\Psi,\varphi\rangle, \ k\in\mathbb{Z}.$ So, if $\Psi$ is a periodic function with period $2L,$  we have 
\[\widehat{\Psi}(k)=\frac 1{2L}\int_{-L}^L \Psi(x)e^{-\frac{ik\pi x}{L}}dx.\]
For $s\in \mathbb{R},$ the Sobolev space of order $s,$ denoted by $H^s_{per}([-L,L])$ is the set of all $f\in \mathcal{P}'$ such that $(1+|k|^{2})^{\frac{s}{2}}\widehat{f}(k)\in l^2(\mathbb{Z}),$ with norm
\[||f||^2_{H^s_{per}}=2L\sum_{k=-\infty}^{\infty}(1+|k|^{2})^s|\widehat{f}(k)|^2.\]
We also note that  $H^s_{per}$ is a Hilbert space with respect to the inner product
\[(f|g)_s = 2L\sum_{n=-\infty}^{\infty}(1+|k|^2)^{s}\widehat{f}(k)\overline{\widehat{g}(k)}\]
In the case $s=0,$  $H^0_{per} $ is a Hilbert space that is isometrically isomorphic to $L^2([-L,L])$ and
\[(f|g)_0 = (f,g) = \int_{-L}^{L} f\overline{g} \ dx.\]
The space $H^0_{per}$ will be denoted by $L^2_{per}$ and its norm will be $\|\cdot\|_{L^2_{per}}.$
Of course $H^s_{per} \subset L^2_{per}$, for any $s \geq 0 $. Moreover, $(H^s_{per})'$, the topological dual of $H^s_{per}$, is isometrically isomorphic to $H^{-s}_{per}$ for all $s \in \mathbb{R}$. The duality is implemented concretely by the pairing
\[\langle f,g\rangle_s = 2L\sum_{k=-\infty}^{\infty}\widehat{f}(k)\overline{\widehat{g}(k)}, \ \ \ for \ \ \ f \in H^{-s}_{per}, \ \  g \in H^s_{per}. \]
Thus, if $f \in L^2_{per}$ and $g \in H^s_{per} $, with $s\geq 0,$ it follows that $\langle f,g\rangle_s = (f,g)$. Additionally, in the particular case $s=\frac 1{2}$ we will denote the pairing $\langle f,g\rangle_s$ simply by $\langle f,g\rangle.$ One of Sobolev's Lemmas in this context states that if $s>\frac{1}{2}$ and
\[C_{per} = \{f: \mathbb{R} \longrightarrow \mathbb{C} \ | \ f \ \ \text{is continuous and periodic with period} \ \  2L \},\]
then $H^{s}_{per}\hookrightarrow C_{per}$.\\

Let $s\in \mathbb{R}.$ The ($L^2$ type) Sobolev space $H^s(\mathbb{R})$ is the collection of all $f\in \mathcal{S}'(\mathcal{R})$ such that $(1+|\xi|^2)^{\frac s{2}}\widehat{f}\in L^2(\mathbb{R},d\xi),$ that is,  $\widehat{f}$ is a measurable function and 
\[
\| f\|_{s}^2=\int_{\mathbb{R}}(1+|\xi|^2)^s|\widehat{f}(\xi)|^2d\xi<\infty.
\]
For more details see Iorio and Iorio \cite{ioriolibro}. Finally, we say that $b\in \widehat{H}^{-1}(\mathbb{R})$ if there exists $V\in L^2(\mathbb{R})$ such that $b=-V'$ and $\|b\|_{\widehat{H}^{-1}}=\|V\|_{L^2}.$
\section{Existence of dnoidal wave solutions}
In this section we  show the existence of a smooth curve of dnoidal wave solutions, with the same fundamental period, for the Zakharov system. In this case, we are interested in solutions for the system (\ref{equaZakha}) in the form given in (\ref{solforms}). Since $u$ is a periodic function (with period $L$), for $c\neq 0$ we suppose that there exists $m\in\mathbb{N}$ such that $L=\frac{4\pi m}{c}.$ Note that for $c=0$ we obtain immediately that $u$ is a $L$-periodic function. Substituting (\ref{solforms})  in (\ref{equaZakha}), we have that $\phi=\phi_{\omega,c}$ and $\psi=\psi_{\omega,c}$ have to satisfy (\ref{systemedo}). Integrating the first equation in (\ref{systemedo}) we obtain
\begin{equation}\label{edo2Zakha}
(c^2-1)\psi'=(\phi^2)'+a_0.
\end{equation}
Using the fact that $\phi^2$ and $\psi$ are periodic we get that $a_0=0.$ Therefore 
\begin{equation}\label{ecuacorr}
(c^2-1)\psi'=(\phi^2)'.
\end{equation}
Integrating (\ref{ecuacorr}), we have that for all $c\neq 1$
\begin{equation}\label{ecuasegint}
\psi=\frac{-\phi^2}{1-c^2}+a_1.
\end{equation}
We assume in our theory that the constant of integration $a_1$ is  zero. Thus, substituting (\ref{ecuasegint}) in the second equation of (\ref{systemedo}) we have that
\begin{equation} \label{ecuaordphi}
\phi''+\left(\omega+\frac{c^2}{4}\right)\phi+\frac{\phi^3}{1-c^2} = 0.
\end{equation}
Now, multiplying (\ref{ecuaordphi}) by $\phi '$ and integrating once, we arrived at
\[\frac{\left(\phi'\right)^2}{2}+\left(\omega+\frac{c^2}{4}\right)\frac{\phi^2}{2}+\frac{\phi^4}{4(1-c^2)}=A_{\phi},\]
where  $A_{\phi}$ is a constant of integration. Then,
\[\left(\phi'\right)^2=\frac{1}{2(1-c^2)}F(\phi),\]
 where $F$ is a polynomial given by
\[F(t)=-t^4+2(1-c^2)\left(-\omega-\frac{c^2}{4}-a_1\right)t^2+4(1-c^2)A_{\phi}.\]
Suppose that $F$ has roots $\pm\eta_1$ and  $\pm\eta_2$ (note that $F$ is even) and without loss of generality that $0<\eta_2<\eta_1$. Thus, we can write
\begin{equation} \label{equa9zak}
\left(\phi'\right)^2=\frac{1}{2(1-c^2)}(\phi^2-\eta^2_2)(\eta^2_1-\phi^2).
\end{equation}
Assume also that $1-c^2>0,$ then the left side of (\ref{equa9zak}) is not negative, therefore we have that
\[\eta^2_2 \leq \phi^2\leq \eta^2_1.\]
Since we are interested in positive solutions, from the last inequality we obtain $\eta_2\leq\phi\leq\eta_1.$  Using (\ref{equa9zak}) we get that the $\eta_{j}$'s satisfy
\[\left \{
\begin{aligned}
-2(1-c^2)\left(\omega+\frac{c^2}{4}\right)&=\eta^2_1 + \eta^2_2\\
4(1-c^2)A_{\phi}&=-\eta^2_1\eta^2_2.\\
\end{aligned}\right.\]
From the last system, we get the restriction  $4\omega+c^2<0.$ \\

Now, define $\varrho(\xi)=\frac{\phi(\xi)}{\eta_1}$, $k^2=\frac{\eta^2_1-\eta^2_2}{\eta^2_1}$ and assume that $\varrho(0)=1$. Thus, we can rewrite the equation (\ref{equa9zak}) as
\begin{equation} \label{equa11zak}
\left(\varrho'\right)^2=\frac{\eta^2_1}{2(1-c^2)}(1-\varrho^2)(\varrho^2+k^2 -1).
\end{equation}
Finally, define $\chi$ through the relation $\varrho^2=1-k^2\sin^2\chi$, with $\chi(0)=0,$ then (\ref{equa11zak}) can be reduce to 
\begin{equation}\label{ecuacorr2}
[\chi']^2=\frac{\eta^2_1}{2(1-c^2)}(1-k^2\sin^2\chi).
\end{equation}
From (\ref{ecuacorr2}) we obtain after some algebra that,
\begin{equation}\label{ecuacorr3}
\int^{\chi(\xi)}_0 \frac{dt}{\sqrt{1-k^2\sin\ t}}= \frac{\eta_1\xi}{\sqrt{2(1-c^2)}}.
\end{equation}
Using the identity (\ref{ecuacorr3}), we obtain from the definition of the Jacobi elliptic functions (see Byrd and Friedman \cite{byrdFriedman}) that
\[\sin(\chi(\xi))=\text{sn}\left(\frac{\eta_1\xi}{\sqrt{2(1-c^2)}} ;k\right).\]
Therefore
\begin{align*}
\varrho(\xi)&=\sqrt{1-k^2\sin^2(\chi(\xi))}=\sqrt{1-k^2 \text{sn}^2\left(\tfrac{\eta_1\xi}{\sqrt{2(1-c^2)}}  ;k\right)} = \text{dn}\left(\tfrac{\eta_1\xi}{\sqrt{2(1-c^2)}} ;k\right),
\end{align*}
where we use the fact that $k^2\text{sn}^2+\text{dn}^2=1.$ Coming back to the variable $\phi$ we obtain that
\begin{equation} \label{equa12zak}
\phi_{\omega,c}(\xi) = \eta_1\text{dn}\left(\tfrac{\eta_1\xi}{\sqrt{2(1-c^2)}} ;k\right)
\end{equation}
and using (\ref{ecuasegint}) we arrive at
\begin{equation} \label{equa13zak}
\psi_{\omega,c}(\xi)=-\frac{\eta^2_1}{1-c^2}\text{dn}^2 \left(\tfrac{\eta_1\xi}{\sqrt{2(1-c^2)}};k\right).
\end{equation}

Now, since dn has fundamental period  $2K$, where  $K = K(k)$ is the complete elliptic integral of the first type (see Byrd and Friedman \cite{byrdFriedman}), we obtain that  $\phi$ e $\psi$ have fundamental period given by
\[T_\psi = T_\phi = \frac{2{\sqrt{2(1-c^2)}}}{\eta_1}K(k).\]
Fix $\omega$ and $c$ such that $1-c^2 >0$  and  $4\omega+c^2<0.$ Additionally, define 
\[\nu=-\left(\omega+\frac{c^2}{4}\right) \ \  \text{and} \ \ \alpha=1-c^2.\] 
Then $\eta^2_1 + \eta^2_2 = 2\nu \alpha$ and consequently $0 < \eta_2 < \sqrt{\nu\alpha} < \eta_1 < \sqrt{2\nu\alpha}.$ We express  $T_\psi$ and $T_{\phi}$ as  functions of the parameter $\eta_2,$
\[T_\psi(\eta_2)= T_\phi(\eta_2)=\frac{2{\sqrt{2\alpha}}}{\sqrt{2\nu\alpha-\eta^2_2 }}K(k(\eta_2)) \ \ \ \ \text{with}\ \ \ \ \ 
k^2(\eta_2)=\frac{2\nu\alpha-2\eta^2_2}{2\nu\alpha-\eta^2_2}.\]
Note that if $ \eta_2\rightarrow 0$, we have that $k(\eta_2)\rightarrow 1^-,$ which implies that $K(k(\eta_2))\rightarrow +\infty $ and consequently $ T_\psi,T_\phi\rightarrow +\infty $. On the other hand, when $\eta_2\rightarrow {\sqrt{\nu\alpha}}$ we get that $k(\eta_2)\rightarrow 0^+$ and then $K(k(\eta_2))\rightarrow \frac{\pi}{2}$. Therefore  $T_\psi,T_\phi \rightarrow \frac{\pi{\sqrt{2}}}{\sqrt{\nu}}.$ Since the function $\eta_2 \in (0,\sqrt{\nu\alpha})\mapsto T_\psi(\eta_2)=T_\phi(\eta_2)$ is strictly decreasing  (we prove this fact later) we obtain
\[T_\phi=T_\psi > \frac{\pi{\sqrt{2}}}{\sqrt{\nu}}.\]  Now, for  $L>0$ and $1-c^2>0$ fixed, chose $\nu > 0$ such that $\sqrt{\nu}> \frac{\pi{\sqrt{2}}}{L}.$
Then, it follows from the analysis given above that there exists a unique $\eta_2=\eta_2(\nu)\in(0,\sqrt{\nu\alpha})$ such that the dnoidal waves $\phi=\phi(\cdot;\eta_1(\nu),\eta_2(\nu))$ and $\psi=\psi(\cdot;\eta_1(\nu),\eta_2(\nu))$ have fundamental period $L=T_{\psi}(\eta_2)=T_{\phi}(\eta_2).$
\begin{remark}
The  formula  (\ref{equa12zak}) and (\ref{equa13zak}) contains, at least formally, the solitary wave solutions for the system  (\ref{equaZakha}) found by Ya Ping in \cite{Wu1}. In fact,  if $\eta_2 \rightarrow 0^+$ we obtain that  $\eta_1 \rightarrow \sqrt{2\nu\alpha}$, $k(\eta_2) \rightarrow 1^-$ and $dn(x,1^-)=\text{sech}(x)$. Consequently 
\[\phi_{c,\omega}(x)=\frac{\sqrt{(-4\omega - c^2)(1-c^2)}}{2} \text{sech} \left(\frac{\sqrt{-4\omega - c^2}}{2} \ x\right)\ \ \text{and} \]
\[\psi_{c,\omega}(x)=\left(2\omega+\frac{c^2}{2}\right) \text{sech}^2 \left(\frac{\sqrt{-4\omega - c^2}}{2} \ x\right).\]
\end{remark}
\begin{theo}\label{TeorDnoidal}
Let $L>0$ and $1-c^2>0$ be arbitrarily  fixed. Consider $\nu_0>\frac{2\pi^2}{L^2}$ and the unique  $\eta_{2,0}=\eta_{2,0}(\nu_0)\in (0,\sqrt{\nu_0\alpha})$ such that $T_{\psi_{\nu_0}}=L=T_{\phi_{\nu_0}}$. Then,

(i) there exist intervals $I(\nu_0)$ and $B(\eta_{2,0})$ around $\nu_0$  and $\eta_{2,0}$ respectively, and a unique smooth function  $\Lambda:I(\nu_0)\longrightarrow B(\eta_{2,0})$ such that $\Lambda(\nu_0)=\eta_{2,0}$ and
\[\frac{2\sqrt{2\alpha}}{\sqrt{2\nu\alpha - \eta_2^2}}K(k) = L,\]
for all $\nu\in I(\nu_0)$, $\eta_2=\Lambda(\nu)$ and
\begin{equation}\label{equaZa13c}
k^2 = k^2(\nu)=\frac{2\nu\alpha-2\eta_2^2}{2\nu\alpha-\eta_2^2}.
\end{equation}
Furthermore, we can chose  $I(\nu_0)=(\frac{2\pi^2}{L^2},+\infty).$

(ii) The dnoidal waves $\psi(\cdot;\eta_1,\eta_2)$ and $\phi(\cdot;\eta_1,\eta_2)$ given by  (\ref{equa12zak}) and (\ref{equa13zak}), and  determined by  $\eta_1=\eta_1(\nu),$  $\eta_2=\eta_2(\nu)=\Lambda(\nu),$ with  $\eta^2_1 + \eta^2_2 = 2\nu\alpha,$ have fundamental period $L$ and satisfy (\ref{ecuasegint}) and (\ref{ecuaordphi}). Furthermore, the map
\[\nu \in I(\eta_0)\longmapsto\left(\psi(\cdot;\eta_1(\nu),\eta_2(\nu)),\phi(\cdot;\eta_1(\nu),\eta_2(\nu))\right)\in H_{per}^n([0,L])\times H_{per}^n([0,L]) \]
is smooth for all integer  $n\geq 1$.

(iii) The map  $\Lambda:I(\nu_0)\rightarrow B(\eta_{2,0})$ is strictly decreasing. Therefore, from (\ref{equaZa13c}), $\nu\mapsto k(\nu)$ is a strictly increasing function.
\end{theo}
\proof
The proof of this theorem follows the same ideas of the Theorem  2.1 in Angulo \cite{angulo4}, we will use the Implicit Function Theorem. For this, consider the open set 
\[\Omega = \left\{(\eta,\nu) \in \mathbb{R}^2: \nu>\frac{2\pi^2}{L^2} \ \text{and} \ \eta \in (0,\sqrt{\nu\alpha}\ )\right\}\] and $\Gamma:\Omega \longrightarrow \mathbb{R}$ defined as
\[\Gamma(\eta,\nu) = \frac{2\sqrt{2\alpha}}{\sqrt{2\nu\alpha-\eta^2}} \ K(k(\eta,\nu))- L,\]
where
\begin{equation}\label{equa13bzak}
k^2(\eta,\nu) =\frac{2\nu\alpha-2\eta^2}{2\nu\alpha-\eta^2}.
\end{equation}
From the hypothesis, we have that $\Gamma(\eta_{2,0},\nu_0) = 0$. We proof that $\frac{d\Gamma}{d\eta}<0$ in $\Omega$. In fact, we use the next relation
\begin{equation}\label{equa14zak}
\frac{dK(k)}{dk}=\frac{E(k)-k'^2K(k)}{kk'^2}\ \ \text{with}\ \ k \in (0,1),
\end{equation}
where $E=E(k)$ is the complete elliptic integral of the second type and  $k'^2=1-k^{2}$ is the complementary modulus. Deriving  (\ref{equa13bzak}) with respect to  $\eta,$ we obtain that 
\begin{equation} \label{equa15zak}
\frac{\partial k}{\partial\eta}= -\frac{2\eta\nu\alpha}{k(2\nu\alpha-\eta^2)^2}.
\end{equation}
Then from (\ref{equa14zak}) and (\ref{equa15zak}) we obtain
\[\frac{\partial\Gamma}{\partial\eta}=\frac{2\eta\sqrt{2\alpha}}{(2\nu\alpha-\eta^2)^{\frac{3}{2}}} \ K(k)-\frac{4\eta\nu\alpha\sqrt{2\alpha}}{(2\nu\alpha-\eta^2)^{\frac{5}{2}}}\left[\frac{E(k) - k'^2K(k)}{k^2k'^2}\right].\]
Thus,
\begin{align*}
\frac{\partial \Gamma}{\partial \eta} < 0 & \Leftrightarrow k^2k'^2(2\nu\alpha-\eta^2)K(k) < 2\nu\alpha E(k) - 2\nu\alpha k'^2K(k)\\
& \Leftrightarrow k'^2(2\nu\alpha-2\eta^2)K(k) + 2\nu\alpha k'^2K(k) <  2\nu\alpha E(k)\\
& \Leftrightarrow \frac{2\nu\alpha k'^2}{(1+k'^2)}K(k)<  \nu\alpha E(k)\Leftrightarrow (1+k'^2)E(k)- 2k'^2K(k) > 0.
\end{align*}
Since the last inequality always holds, we obtain that  $\frac{\partial \Gamma}{\partial \eta}<0$. By the Implicit Function Theorem we have that there exists an interval  $I(\nu_0)$ around  $(\nu_0)$, an interval  $B(\eta_{2,0})$ around $\eta_{2,0}$  and a smooth function $\Lambda:I(\nu_0)\longrightarrow B(\eta_{2,0})$ such that  $\Lambda (\nu_0) = \eta_{2,0}$ and
\[\Gamma(\Lambda(\nu),\nu) = 0,\ \ \ \ \ \ \forall \nu \in I(\nu_0).\]
Additionally, since  $\nu_0$ was chosen arbitrarily in  $I=\left(\frac{2\pi^2}{L^2},+\infty\right)$ and from the uniqueness of  $\Lambda,$ we extend $\Lambda$ to  $I.$  The part $(ii)$ is immediate, using the smoothness of the function involved.\\

Now, we prove that  $\Lambda$ is an strictly decreasing function. For this note that  $\Gamma(\Lambda(\nu),\nu)=L$ for all  $\nu\in I(\nu_0)$ then, using again the Implicit Function Theorem we get that 
\[\Lambda '(\nu)=-\frac{\partial\Gamma/{\partial\nu}}{\partial\Gamma/\partial\eta}.\]
Since $\frac{\partial\Gamma}{\partial\eta}<0,$ we just have to prove that  $\frac{\partial\Gamma}{\partial\nu}<0$ in  $I(\nu_0)$. In fact, since 
\[\frac{\partial\Gamma}{\partial\nu}=\frac{2\alpha\sqrt{2\alpha}}{(2\nu\alpha-\eta^2)^{3/2}}\left[ -K+\frac{dK}{dk}\frac{\eta^2}{k(2\nu\alpha-\eta^2)} \right]\]
and  $\eta^2=(2\nu\alpha-\eta^2)k'^2,$ we obtain
\begin{align*}
\frac{\partial\Gamma}{\partial\nu}<0&\Leftrightarrow \frac{\eta^2}{\sqrt{2\nu\alpha-\eta^2}\sqrt{2\nu\alpha-2\eta^2}}\frac{dK}{dk}<K\Leftrightarrow k'^2\frac{dK}{dk}-kK<0.\\
\end{align*}
From  (\ref{equa14zak}), we arrived at
\[ \frac{\partial\Gamma}{\partial\nu}<0\Leftrightarrow\frac{E-k'^2K}{k}-kK<0\Leftrightarrow E<K.\]
Since the last inequality always holds for any $k\in(0,1)$  (see  Byrd and Friedman \cite{byrdFriedman}), we obtain the desired result.\\

Finally, deriving $k$ with respect to  $\nu,$ we obtain
\[\frac{dk}{d\nu}=\frac{\alpha\eta(\eta-2\eta'\nu)}{k(2\nu\alpha-\eta^2)^2}>0,\]
which proves that  $\nu \mapsto k(\nu)$ is strictly increasing function, this finishes the proof of the theorem.
\endproof

The next result will be used in the prove of the stability of the  dnoidal waves solutions.
\begin{coro}\label{coroDnoidal}
Let $L>0$ and $c$ be arbitrarily fixed with $1-c^2 >0.$  Consider the smooth curve  of dnoidal waves $\nu \in \left(\frac{2\pi^2}{L^2},+\infty\right) \longmapsto \phi_{\nu}(\cdot ; \eta_1(\nu),\eta_2(\nu))$ determined by Theorem \ref{TeorDnoidal}. Then
\[\frac{d}{d\nu}\int_0^L \phi^2_\nu(\xi)d\xi>0.\]
\end{coro}
\proof
Using the facts that $\eta_1L=2\sqrt{2(1-c^2)}K(k)$ and $\int_0^L dn^2 (y) dy = E(k)$ (see Byrd and Friedman \cite{byrdFriedman}) we get that
\[\int_0^L \phi^2_{\nu}(\xi) d\xi = 2\eta_1\sqrt{2(1-c^2)} \int_0^K dn^2 (y;k) dy = \frac{8(1-c^2)}{L}K(k)E(k).\]
Since  $k\mapsto K(k)E(k)$ and $\nu \mapsto k(\nu)$ are strictly increasing functions we obtain
\[\frac{d}{d\nu}\int_0^L \phi^2_{\nu}(\xi) d\xi =\frac{8(1-c^2)}{L} \frac{d}{dk }\left[K(k)E(k)\right]\frac{dk}{d\nu}>0\]
This finishes the proof of the corollary.
\endproof
\section{Spectral Analysis}
In this part of the paper, we study some spectral properties of various operators  which will be necessary to obtain our result of stability. First, note that the system  (\ref{equaZakha}) can be rewritten as 
\begin{equation}\label{ZakNovo}
\left \{
\begin{aligned}
v_t &=-V_x, \ \int_0^L V(x,t)dx  = 0 \\
 V_t &= -(v+|u|^2)_x \\
iu_t &+ u_{xx} =uv
\end{aligned} \right.
\end{equation}
Therefore, we have the Hamiltonian structure $\frac{\partial U}{\partial t}=JE'(U)$ where $U=(v,V,u)^t,$ $J$ is the linear skew-symmetric  operator given by 
\[J=\left(
\begin{array}{crrc}
 0 &-\frac{d}{dx} &0 \\
 -\frac{d}{dx} &0 &0 \\
 0 &0 &-\frac{i}{2}
\end{array}
\right)\]
and $E$ is the energy functional define as
\begin{equation} \label{equa2.4Zak}
E(v,V,u)=\frac{1}{2}\int_0^L 2|u_x|^2+v^2+V^2+2v|u|^2 \ dx.
\end{equation}
We also use the functionals $Q_1$ and $Q_2$ defined as
\begin{equation}\label{quantCons2e3}
Q_1(v,V,u)=\int_0^L uV + \text{Im}(u_x\overline{u})\ dx \ \ \ \ \text{and}\ \ \ \ Q_2(v,V,u)= \int_0^L |u|^2 dx.
\end{equation}
A standard analysis proves that $E,$ $Q_1$  and $Q_2$ are conserved quantities of the system  (\ref{ZakNovo}), i.e.,
\[E(v(t),V(t),u(t))= E(v(0),V(0),u(0)), \ \ \ Q_1(v(t),V(t),u(t))= Q_1(v(0),V(0),u(0))\]
\[\text{and}\ \  Q_2(v(t),V(t),u(t))= Q_2(v(0),V(0),u(0)) \]
for all $t\in [-T,T],$ where $T$ is the maximal time of existence of  solutions.\\

Now, suppose that $V(x,t)=\varphi_{\omega,c}(x-ct),$ with $\varphi_{\omega,c}:\mathbb{R}\rightarrow\mathbb{R}$ a smooth  $L-$periodic function, is solution of $v_t=-V_{x}$, then
\begin{equation}\label{edo3Zakh}
c\psi_{\omega,c}'=\varphi'_{\omega,c}.
\end{equation}
Therefore, $c\psi= \varphi+d_0$, where $d_0$ is a constant of integration. Since we are interested in  $\varphi$ with zero mean, we obtain that $d_0=\frac{c}{L}\int_0^L \psi dx$. Using tha fact that $\int_0^K \text{dn}^2(x,k) dx=E(k)$ we get that $d_0(k)=-\frac{c\eta_1^2}{1-c^2}\frac{E(k)}{K(k)}$. Therefore
\begin{equation}\label{varphidnoidal}
\varphi(\xi)=-\frac{c\eta_1^2}{1-c^2}\left[\text{dn}^2\left(\frac{\eta_1\xi}{\sqrt{2(1-c^2)}};k\right)-\frac{E(k)}{K(k)}\right].
\end{equation}

It is worth to note that if $\eta_2\rightarrow 0^+,$ then $\eta_1\rightarrow\sqrt{2\alpha\nu}$ and therefore $k\rightarrow 1^-.$ Since dn$(u,1^-)=\text{sech}(u)$, $E(1)=\frac{\pi}{2}$ and $K(1)= +\infty$ we arrive at
\[\varphi(\xi)=c\left(2\omega+\frac{c^2}{2}\right)\text{sech}^2\left(\frac{\sqrt{-4\omega-c^2}}{2}\xi \right),\]
which is the solitary wave solution for (\ref{edo3Zakh}).\\

Now, using the Theorem \ref{TeorDnoidal} we have that there exist periodic traveling waves for  (\ref{ZakNovo}) given by
\[\left(\psi_{\omega,c}(x-ct),\varphi_{\omega,c}(x-ct), e^{- i\omega t} e^{i\frac{c}{2}(x-ct)}\phi_{\omega,c}(x-ct)\right),\]
where
\begin{equation}\label{solZakhpsi}
\psi_{\omega,c}(\xi)=\frac{-\eta_{1}^2}{1-c^2}\text{dn}^2\left(\frac{\eta_1\xi}{\sqrt{2(1-c^2)}};k \right),
\
\phi_{\omega,c}(\xi)=\eta_1\text{dn}\left(\frac{\eta_1\xi}{\sqrt{2(1-c^2)}};k\right)
\end{equation}
\begin{equation}\label{solZakhphi}
\text{and}\ \ \ \ \varphi_{\omega,c}(\xi)= -\frac {c\eta_1^2}{1-c^2}\left[\text{dn}^2\left(\frac{\eta_1\xi}{\sqrt{2(1-c^2)}};k\right)-\frac{E(k)}{K(k)}\right].
\end{equation}
The next operators will be useful in the proof of the stability of the dnoidal wave solutions:
\begin{equation}\label{operaDnoidal}
\mathcal{L}_{3}=-\frac{d^2}{{dx}^2}-\left(\omega+\frac{c^2}{4}\right)+3\psi\ \ \ \text{and}\ \ \ \
\mathcal{L}_{4}=-\frac{d^2}{{dx}^2}-\left(\omega+\frac{c^2}{4}\right)+\psi.
\end{equation}
We will study the spectral properties of the operators $\mathcal{L}_{i},\ i=3,4$. Recall that  $\sigma(\mathcal{L}_{i})=\sigma_{ess}(\mathcal{L}_{i})\cup\sigma_{disc}(\mathcal{L}_{i})$ where $\sigma_{ess}(\mathcal{L}_{i})$ and $\sigma_{disc}(\mathcal{L}_{i})$  denote, respectively, the  essential spectrum and the point spectrum of $\mathcal{L}_{i}$ (see Reed and Simon \cite{ReedSimon1}). Write
\[\mathcal{L}_{3}=\left(-\frac{d^2}{{dx}^2}-\left(\omega+\frac{c^2}{4}\right)\right)+3\psi=:\mathcal{L}+M_1,\]
\[\mathcal{L}_{4}=\left(-\frac{d^2}{{dx}^2}-\left(\omega+\frac{c^2}{4}\right)\right)+\psi=:\mathcal{L}+M_2,\]
where $\mathcal{L}=-\frac{d^2}{dx^2}-(\omega+\frac{c^2}{4}).$ Since  $M_1$  and $M_2$ are relatively compact with respect to $\mathcal{L}$, it follows from the Weyl's Essential Spectrum Theorem (see Reed and Simon  \cite{ReedSimon1}) that  $\sigma_{ess}(\mathcal{L}_i)=\sigma_{ess}(\mathcal{L})=\emptyset,$ with $i= 3,4.$ Thus $\sigma(\mathcal{L}_i)=\sigma_{disc}(\mathcal{L}_i),$ para $i=3,4.$ Therefore we have to analyze the periodic eigenvalue problem on $[0,L]$
\begin{equation}\label{probperiogeral}
\left \{
\begin{aligned}
\mathcal{L}_{i}\chi&=\lambda\chi\\
\chi(0)&=\chi(L),\ \chi'(0)=\chi'(L). \\
\end{aligned}\right.
\end{equation}
The problem (\ref{probperiogeral}) determines that the spectrum of  $\mathcal{L}_i$ is a countable set of eigenvalues  $\{\lambda_n: n=0,1,2,3,...\}$ with 
\[\lambda_0\leq \lambda_1\leq \lambda_2\leq \lambda_3\leq  \cdots, \]
where the double eigenvalues are counted twice and $\lambda\rightarrow+\infty$ when $n\rightarrow\infty.$ We denote by  $\chi_n$ the eigenfunctions  associated to  the eigenvalue $\lambda_n.$ It is clear from the conditions $\chi(0)=\chi(L), \chi'(0)=\chi'(L)$ that $\chi_n$ can be extended to all $(-\infty,+\infty)$ as a continuous differentiable function with period $L.$ We know from the Floquet theory that the periodic eigenvalue problem  (\ref{probperiogeral}) is related to the study of the next semi-periodic eigenvalue problem consider in $[0,L]$
\[\left \{
\begin{aligned}
 \mathcal{L}_{i}\eta&= \mu\eta\\
\eta(0)&=-\eta(L),\ \eta'(0)=-\eta'(L), \\
\end{aligned} \right.\]
which also is a self-adjoint problem and therefore determines a sequence of eigenvalues $\{\mu_n: n=0,1,2,3 ...\}$ with
\[\mu_0\leq \mu_1\leq \mu_2\leq \mu_3\leq  \cdots, \]
where the double eigenvalues are counted twice and $\mu_n\rightarrow+\infty$ when $n\rightarrow\infty.$
We denote by  $\eta_n$ the eigenfunction associated to the eigenvalue $\mu_n.$
\begin{theo}
Let $\phi_{\nu}=\phi$ and $\psi_{\nu} = \psi$ the dnoidal waves given by  Theorem \ref{TeorDnoidal}. Then,

(i) the operator $\mathcal{L}_3$  in (\ref{operaDnoidal}) defined in $L_{per}^2([0,L])$ with domain  $H_{per}^2([0,L])$ has its fist three eigenvalues simple, where zero is the second one with associated eigenfunction $\phi '$. Furthermore, the rest of the spectrum is constitute by a discrete set of eigenvalues which are double.

(ii) The operator $\mathcal{L}_4$ in (\ref{operaDnoidal}) defined in  $L_{per}^2([0,L])$ with  domain $H_{per}^2([0,L])$ has zero as its first eigenvalue which is simple with associated eigenfunction $\phi$. Furthermore, the rest of the spectrum is constitute by a discrete set of eigenvalues.
\end{theo}
\proof
$(i)$ The proof is based on the Floquet Theory (see Eastham  \cite{eastham},  Mangnus and Winkler \cite{magnus}). Deriving (\ref{ecuaordphi}) and using (\ref{ecuasegint}) we have that  $\mathcal{L}_3\phi '=0$. Then zero is an eigenvalue of $\mathcal{L}_3$ with associated eigenfunction  $\phi '.$ Since $\phi '$ has exactly two zeros  on $[0,L)$, we get that zero is either the second or the third eigenvalue of $\mathcal{L}_3$. We will prove that zero is in fact the second one. For this we have to study the periodic problem
\begin{equation} \label{equa2.12Zak}
\left \{
\begin{aligned}
\mathcal{L}_{3}\chi &=\lambda\chi\\
\chi(0)&=\chi(L),\ \chi'(0)=\chi'(L).\\
\end{aligned} \right.
\end{equation}
Let $\Lambda(x)=\chi(\eta x)$ where  $\eta=\frac{\sqrt{2\alpha}}{\eta_1}$. Then, from the explicit form of $\psi$ and the relation  $k^2 \text{sn}^2+\text{dn}^2=1$,we have that the problem  (\ref{equa2.12Zak}) is equivalent to 
\begin{equation} \label{equa2.13Zak}
\left \{
\begin{aligned}
\Lambda ''+&[\rho-6k^2\text{sn}^2(x; k)]\Lambda=0\\
\Lambda(0)&=\Lambda(2K),\ \Lambda '(0)=\Lambda '(2K),  \\
\end{aligned}\right.
\end{equation}
where
\[\rho=\frac{2\alpha}{\eta_{1}^2}\left(\frac{\lambda}{2}+\omega+\frac{c^2}{4}+\frac{3\eta_{1}^2}{\alpha}\right).\]
The second order equation given in (\ref{equa2.13Zak}) is called the Jacobian form of the Lam\' e equation. It is well known that such equation determines  the existence of exactly three intervals of instability (see Theorem $7.8$ in Mangnus and Winkler \cite{magnus}). We will show that this intervals are the first three. First, observe that $\rho_1= 4+k^2$ and $\Lambda_1(x)=\text{cn}(x;k)\text{sn}(x;k)$ satisfy the problem (\ref{equa2.13Zak}). Furthermore, following Ince \cite{ince} we have that the functions
\[\Lambda_0(x)=1-(1+k^2-\sqrt{1+k^2+k^4})\text{sn}^2(x;k),\]
\[\Lambda_2(x)=1-(1+k^2+\sqrt{1+k^2+k^4})\text{sn}^2(x;k),\]
which have period $2K,$ are the eigenfunctions of  (\ref{equa2.13Zak}) with eigenvalues given by
\[\rho_0= 2\left(1+k^2 -\sqrt{1+k^2+k^4}\right)\ \ \ \text{and}\ \ \  \rho_2= 2\left(1+k^2 -\sqrt{1+k^2+k^4}\right).\]
Since  $\Lambda_0$ does not have zeros in $[0,2K]$, it follows that  $\rho _0$ is the first eigenvalue of  (\ref{equa2.13Zak}). Furthermore, since $\Lambda_2$ has two zeros in $[0,2K)$ and $\rho_1<\rho_2,$ we have that  $\rho_1$ is the second eigenvalue of  (\ref{equa2.13Zak}) and $\rho_2$ is the third. We also have that $\rho_0,\rho_1$ and $\rho_2$ are simple. Now, since the eigenvalues of  (\ref{equa2.12Zak}) and (\ref{equa2.13Zak}) are related as
\[\lambda=\frac{\eta_1}{\alpha}(\rho-6)+2\nu,\]
we can see  $\lambda$ as a function of  $\rho$, which is increasing. Since $k^2-2=\frac{2\alpha\nu}{\eta_{1}^2}$, we have that  $\lambda(\rho_1) = 0 = \lambda_1$ and since $\lambda_0 < \lambda_1 < \lambda_2$, we obtain that 
\[\lambda_0<0=\lambda_1<\lambda_2.\]
This finishes the proof of the part  $(i)$.\\

$(ii)$ Using (\ref{ecuaordphi}) and (\ref{ecuasegint}) we have that  $\mathcal{L}_4\phi=0$. Thus,  zero is an eigenvalue of  $\mathcal{L}_4$ with associated eigenfunction $\phi$. Since  $\phi$ does not have zeros in  $[0,L]$ we obtain that zero is the first eigenvalue of  $\mathcal{L}_4$ and it is simple.
\endproof
\section{Nonlinear Stability for the Dnoidal Wave Solutions}
In this section we study the nonlinear stability properties of the periodic traveling wave solution $\Phi(\xi)=(\psi(\xi),\varphi(\xi),\widetilde{\phi}(\xi))$ where $\psi$, $\varphi$ and $\phi$ are given by  (\ref{solZakhpsi}),  (\ref{solZakhphi}), $\widetilde{\phi}(\xi)=e^{i\frac{c}{2}\xi}\phi(\xi)$ and $1-c^2>0$. First, we define the type of stability in which we are interested: Let  $ X:= L^2_{per}([0,L]) \times \widetilde{L}^2_{per}([0,L]) \times H^1_{per}([0,L]),$ where
\[\widetilde{L}^2_{per}([0,L])=\left\{f \in L^2_{per}([0,L]): \int_0^L f(x) dx=0\right\}.\]
Initially, observe that the system (\ref{ZakNovo}) has two basic symmetries: translations and rotations. This means that if  $(v(x,t),V(x,t),u(x,t))$ is a solution of (\ref{ZakNovo}), then  the pair of functions 
\[(v(x+y),V(x+y),u(x+y))\ \ \ \ \ \text{and}\ \ \ \  (v(x,t),V(x,t),e^{-is}u(x,t))\]
are also solutions, for any real constants  $y$ and $s.$ So, our notion of stability will be modulus these symmetries. More precisely,  
\begin{defi}
We say that the orbit generated by $\Phi(\xi)$, namely
\[\mathcal{O}_{\Phi}=\left\{\left(\psi(\cdot + y)),\varphi(\cdot+y)), e^{i\theta}\widetilde{\phi}(\cdot+y)\right):(\theta,y)\in [0,2\pi)\times\mathbb{R}\right\}\]
is stable in $X$ by the flow generated by the system  (\ref{ZakNovo}), if for all  $\epsilon>0$, there exists  $\delta>0$ such that for any $(v_0,V_0, u_0)\in X$ satisfying
\[\|v_0-\psi\|_{L^2_{per}}<\delta,\ \ \|V_0-\varphi\|_{L^2_{per}}<\delta \ \ \text{and}\ \ \|u_0-\widetilde{\phi}\|_{H^1_{per}}<\delta,\]
we have that the solution $(v, V,u)$ of the system  (\ref{ZakNovo}) with $(v(0),V(0),u(0))=(v_0, V_0,u_0)$, satisfies
\[(v,V,u) \in C(\mathbb{R}; L^2_{per}([0,L]))\times C(\mathbb{R}; \widetilde{L}^2_{per}([0,L])) \times C(\mathbb{R};H^1_{per}([0,L])),\]
\begin{equation}\label{desiImp1}
\inf_{y \in \mathbb{R}} \|v(\cdot+y,t)-\psi\|_{L^2_{per}}<\epsilon, \ \ \ \inf_{y \in \mathbb{R}} \|V(\cdot+y,t)-\varphi\|_{L^2_{per}}<\epsilon
\end{equation}
\begin{equation}\label{desiImp2}
\text{and}\ \ \ \inf_{\theta \in [0,2\pi),y \in \mathbb{R}} \|e^{i\theta}u(\cdot+y,t)-\widetilde{\phi}\|_{H^1_{per}}<\epsilon.
\end{equation}
Otherwise, we say that  $\Phi$ é $X$-unstable.   
\end{defi}
Next, we  present our result of stability for the dnoidal waves.
\begin{theo}\label{teoStabDnoidal}
Let $L>0$ and $1-c^2>0$ be fixed numbers. Consider the smooth curve of periodic traveling wave solutions for  the system (\ref{ZakNovo}), $\nu\mapsto(\psi_{\nu},\varphi_{\nu},\phi_{\nu}),$ determined by the Theorem \ref{TeorDnoidal} and (\ref{varphidnoidal}). Then, for  $\nu>\frac{2\pi^2}{L^2}$ the orbit generated by  $\Phi_{\nu}(x,t)=\left(\psi_{\nu}(x),\varphi_{\nu}(x),\widetilde{\phi}_{\nu}(x)\right)$ is stable in $X$ by the periodic flow generated by the system  (\ref{ZakNovo}), if the initial datum $(v_0,V_0,u_0)$ satisfies
\[\int_0^L v_0(x)dx\leq\int_0^L\psi(x)dx.\]
\end{theo}
\proof Consider  $(\psi_{\nu},\varphi_{\nu},\widetilde{\phi}_{\nu})$ the solution of  (\ref{ZakNovo}) given by  Theorem \ref{TeorDnoidal}. For $(v_0,V_0,u_0) \in L^2([0,L]) \times \widetilde{L}_{per}^2([0,L]) \times H_{per}^1([0,L])$ and $(v,V,u)$ the global solution for  (\ref{ZakNovo}) corresponding to this initial data, we define for  $t\geq 0$ and $ \nu > \frac{2\pi^2}{L^2}$
\[\Omega_t(y,\theta)=\|e^{i\theta}(T_cu)'(\cdot + y,t)-\phi_{\nu}'\|_{L_{per}^2}^2+\nu\|e^{i\theta}(T_cu)(\cdot + y,t)-\phi_{\nu}\|_{L_{per}^2}^2,\]
where we denote by $T_c$ the bounded linear operator define as
\[(T_cu)(x,t)=e^{-ic(x-ct)/2}u(x,t).\]
Then, the deviation of the solution $u(t)$ from the orbit generated by $\Phi$ is measure by
\begin{equation}\label{equa3}
\rho_{\nu}(u(\cdot,t),\phi_{\nu})^2:=\inf\left\{\Omega_t(y,\theta): (y,\theta)\in[0,L]\times[0,2\pi]\right\}.
\end{equation}
Therefore, from (\ref{equa3}) we have that for each $t$ the $\inf\Omega_t(y,\theta)$ is attained in $(\theta,y)=(\theta(t),y(t)).$
Consider the perturbation of the periodic wave $(\psi,\varphi,\widetilde{\phi})$
\begin{equation} \label{equa4}
\left \{
\begin{aligned}
 \xi(x,t)&=e^{i\theta}(T_cu)(x+y,t)-\phi_{\nu}(x) \\
 \eta(x,t)&= V(x+y,t)-\varphi_{\nu}(x)\\
 \gamma(x,t)&= v(x+y,t)-\psi_{\nu}(x).
\end{aligned} \right.
\end{equation}
By the property of minimum of  $(\theta,y)=(\theta(t),y(t))$, we obtain from  (\ref{equa4}) that $p(x,t)=\text{Re}(\xi(x,t))$ and $q(x,t)=\text{Im}(\xi(x,t))$ satisfy the compatibility relations
\begin{equation} \label{equa5}
\left \{
\begin{aligned}
 \int_0^L q(x,t)\phi_{\nu}(x)\psi_{\nu}(x) dx&= 0 \\
\int_0^L p(x,t)(\phi_{\nu}(x)\psi_{\nu}(x))' dx &= 0. \\
\end{aligned} \right.
\end{equation}\\
Now, consider the continuous functional  $\mathcal{B}$ defined in  $X$ as
\[\mathcal{B}(v,V,u) := E (v,V,u)-cQ_1(v,V,u)-\omega Q_2(v,V,u),\]
where $E$, $Q_1$ and $Q_2$ were  defined in (\ref{equa2.4Zak}) and (\ref{quantCons2e3}). Then, from (\ref{equa4}) and (\ref{equa5}), we get 
\begin{align*}
\Delta\mathcal{B}:=& \ \mathcal{B}(v(t),V(t),u(t))-\mathcal{B}(\psi,\varphi,\widetilde{\phi})\\
=&\left(\mathcal{L}_3p,p\right)+\left(\mathcal{L}_4q,q\right) +\frac{1}{2}\int_0^L \gamma^2+2\gamma(p^2+q^2)-4\psi p^2 + 4\gamma p \phi \ dx \\&  + \frac{1}{2} \int_0^L 2\gamma\psi +\eta^2+ 2\eta\varphi + 2\gamma\phi^2 - 2c\gamma\eta -2c\gamma\varphi -2c\psi\eta \ dx
\end{align*}
where
\[\mathcal{L}_3=-\frac{d^2}{dx^2}-\left(\omega+\frac{c^2}{4}\right)+3\psi \ \ \  \text{and} \ \ \  \mathcal{L}_4=-\frac{d^2}{dx^2}-\left(\omega+\frac{c^2}{4}\right)+\psi.\]
Using the facts that $c\psi-\varphi=d_0$ and $\int_0^L\eta dx=0,$ we obtain 
\begin{align*}
\Delta\mathcal{B}(t)&=\left(\mathcal{L}_3p,p\right)+\left(\mathcal{L}_4q,q\right) + \frac{1}{2}\int_0^L \left[\sqrt{1-c^2}\gamma+\frac{2\phi p}{\sqrt{1-c^2}}+\frac{p^2+q^2}{\sqrt{1-c^2}}\right]^2 dx \\
&+ \frac{1}{2}\int_0^L(c\gamma-\eta)^2dx - \int_0^L \frac{4\phi p(p^2+q^2)}{1-c^2}+\frac{(p^2+q^2)^2}{1-c^2} dx + \int_0^L(c\gamma-\eta)(c\psi-\varphi)dx\\
& =\left(\mathcal{L}_3p,p\right)+\left(\mathcal{L}_4q,q\right) + \frac{1}{2}\int_0^L \left[\sqrt{1-c^2}\gamma+\frac{2\phi p}{\sqrt{1-c^2}}+\frac{p^2+q^2}{\sqrt{1-c^2}}\right]^2 dx \\
&+\frac{1}{2}\int_0^L(c\gamma-\eta)^2 dx-\int \frac{4\phi p(p^2+q^2)}{1-c^2}+\frac{(p^2+q^2)^2}{1-c^2}dx + cd_0\int_0^L\gamma dx.
\end{align*}
Since $cd_0\leq 0,$  $\int_0^Lv_0dx\leq\int_0^L\psi(x)dx$ and $\int v(t,x)dx=\int v_0(x)dx,$ we have that  $cd_0\int_0^L\gamma dx\geq 0.$ Therefore
\begin{align}\label{equa55}
\notag\Delta\mathcal{B}(t)\geq &\left(\mathcal{L}_3p,p\right)+\left(\mathcal{L}_4q,q\right) + \frac{1}{2}\int_0^L \left[\sqrt{1-c^2}\gamma+\frac{2\phi p}{\sqrt{1-c^2}}+\frac{p^2+q^2}{\sqrt{1-c^2}}\right]^2dx \\&+\frac{1}{2}\int_0^L(c\gamma-\eta)^2dx-C_1\|\xi\|_{H^1_{per}}^3-C_2\|\xi\|_{H^1_{per}}^4,
\end{align}
with $C_i>0$, $i=1,2.$\\

The estimates for  $\left(\mathcal{L}_3p,p\right)$ and $\left(\mathcal{L}_4q,q\right)$ will be obtain from the next theorems.
\begin{theo}\label{teoEstForCua1}
Let $1-c^2>0$ and $\nu>\frac{2\pi^2}{L^2}$ fixed numbers. Consider $\phi_{\nu}$ the  dnoidal wave given by Theorem \ref{TeorDnoidal}. Then
\begin{itemize}
\item[(a)]  $\inf\{\left(\mathcal{L}_3f,f\right):\|f\|=1\ \text{and}\ \left(f,\phi_{\nu}\right)=0\}=:\alpha_0=0$
\item[(b)]  $\inf\{\left(\mathcal{L}_3f,f\right):\|f\|=1,\ \left( f,\phi_{\nu}\right)=0\ \text{and}\  \left(f,(\phi_{\nu}\psi_{\nu})'\right)=0\}=:\alpha>0.$
\end{itemize}
\end{theo}
\proof
$(a)$ Since  $\mathcal{L}_3\left(\frac{d}{dx}\phi_{\nu}\right)=0$ and $\left(\frac{d}{dx}\phi_{\nu},\phi_{\nu}\right)=0$, then  $\alpha_0 \leq 0$. We prove that  $\alpha_0 \geq 0$ using Lemma E.1 in Weinstein \cite{weinstein2} (which works on the periodic case). We first show that the infimum is attained. In fact, since $\phi_{\nu}$ is bounded we have that  $\alpha$ is finite, thus there exists  $\{f_j\}\subset H^1_{per}([0,L])$ with $\|f_j\|=1$, $\left( f_j,\phi_{\nu}\right)=0$ and $\lim_{j\rightarrow \infty}\left(\mathcal{L}_3 f_j,f_j\right)=\alpha_0$. Since $\{f_j\}$ is bounded in $H^1_{per}([0,L])$ there exists a subsequence of  $\{f_j\}$, that we denote again $f_j,$ such that  $f_j \rightharpoonup g$ weakly in $H^1_{per}([0,L]),$ then $f_j \rightarrow g$ in $L^2_{per}([0,L])$. Therefore $(g,\phi_{\nu}) =0$ and $(\phi_c f_j,f_j)\rightarrow (\phi g,g)$ when $j \rightarrow + \infty$. So $g \neq 0$ and $\|g'\|_{L^2_{per}}\leq \liminf \|f_j'\|_{L^2_{per}}$.\\

Now, define $f=g/\|g\|_{L^2_{per}}$, then $(f,\phi_{\nu})=0$, $\|f\|_{L^2_{per}}=1$ and
\[\alpha_0 \leq (\mathcal{L}_3 f,f) \leq\frac{\alpha_0}{\|f\|^2_{L^2_{per}}} =\alpha_0.\]
Therefore the infimum is attained. We show now that  $\alpha_0 \geq 0$. In fact, $\mathcal{L}_3$ has the spectral properties required to use Lemma E.1, we need to find  $\chi$ such that $\mathcal{L}_3\chi=\phi_{\nu}$ and $(\chi,\phi_{\nu})\leq 0$. From Theorem \ref{TeorDnoidal} we have that  $\nu\in\left(\frac{2\pi^2}{L^2},+\infty\right)\longmapsto\phi_{\nu}\in H^1_{per}([0,L])$ is of class $C^1$, then differentiating (\ref{ecuaordphi}) with respect to $\nu$ we obtain that  $\chi=-\frac{d}{d\nu}\phi_{\nu}$ satisfies $\mathcal{L}_3\chi =\phi_{\nu}$. Using the  Corollary \ref{coroDnoidal} we obtain that
\[(\chi,\phi_{\nu})=-\frac{1}{2}\frac{d}{d\nu}\int_0^L \phi^2_\nu(\xi)\ d \xi<0\]
Therefore $(\chi,\phi_\nu)<0,$ which proves that  $\alpha_0 \geq 0$. This finishes the proof of part $(a)$.\\

$(b)$ Using the part $(a)$, we have that  $\alpha \geq 0$. Suppose that $\alpha =0$. Using a similar argument as in part  $(a)$ we obtain that there exists  $f \in H^1_{per}([0,L])$ such that  $\|f\|_{L^2_{per}}=1$ and $(f,\phi_\nu)=\left(f,(\phi_\nu\psi)'\right)=0$. Then, from the theory of  Lagrange Multipliers, there exists  $\lambda,\ \theta$ and $\delta$ such that
\[\mathcal{L}_3 f =\lambda f+\theta \phi_\nu + \delta(\phi_\nu\psi_\nu)'.\]
Since $\left(\mathcal{L}_3 f,f\right)=0$, we obtain that  $\lambda=0$. From the fact that $\mathcal{L}_3\phi'_\nu=0$ we have
\[0=\delta \int_0^L \phi_{\nu}'(\phi_\nu\psi_\nu)'  d\xi = -\frac{3\delta}{1-c^2}\int_0^L (\phi_{\nu}')^2 \phi^2_{\nu}\ d\xi.\]
The last inequality implies  $\delta=0,$ thus  $\mathcal{L}_3f=\theta\phi_\nu$. Consider $\chi=-\frac{d}{d\nu}\phi_\nu$, then we get that  $\mathcal{L}_3(f-\theta\chi)=0$, thus
\[0=(f-\theta\chi,\phi_\nu)=-\theta(\chi,\phi_\nu).\]
Therefore $\theta=0$ and consequently there exists $s \in \mathbb{R}\setminus\{0\}$ such that  $f = s\phi_\nu ',$ which is absurd. Therefore $\alpha>0$, which finishes the proof of the theorem.
\endproof
\begin{theo}\label{teoEstForCua2}
Let $1-c^2>0$ and $\nu>\frac{2\pi^2}{L^2}$ be fixed numbers. Consider $\phi_\nu$ and $\psi_\nu$ the dnoidal waves given by Theorem \ref{TeorDnoidal}. Then,
\[\inf \{\left(\mathcal{L}_4 f,f\right):\|f\|_{L^2_{per}}=1\ \ \text{and} \ \ \left(f,\phi_\nu\psi_\nu\right)=0\}=:\beta>0\]
\end{theo}
\proof
From the spectral properties of $\mathcal{L}_4$ is clear that $\mathcal{L}_4$ is a nonnegative operator, therefore $\beta\geq 0$. Suppose that $\beta=0 $. Then, following the same ideas of the proof of Theorem \ref{teoEstForCua1}, we have that the infimum is attained on a admissible function $g\neq 0$ and there exists $(\lambda,\theta)\in\mathbb{R}^2$ such that 
\[\mathcal{L}_4g =\lambda g + \theta \phi_\nu\psi_\nu.\]
Since $\left(g, \phi_\nu\psi_\nu\right) =0$, then $\lambda =0$. Furthermore,
\[0=(\mathcal{L}_4\phi,g)=\theta\int_0^L\phi^2_\nu\psi_\nu\ d\xi,\]
which implies $\theta =0$. Since zero is a simple eigenvalue of  $\mathcal{L}_4$, there exists  $s\in\mathbb{R}\setminus\{0\}$ such that  $g= s\phi$, which is absurd. This finishes the proof of the theorem.
\endproof

Our goal is to estimate the terms  $\left(\mathcal{L}_3 p,p\right)$ and $\left(\mathcal{L}_4 q,q\right)$, where  $p$ and $q$ satisfy (\ref{equa5}). Using  Theorem \ref{teoEstForCua2} and definition  of $\mathcal{L}_4$, we have that there exists  $C_0>0$ such that 
\begin{equation}\label{equa8}
\left(\mathcal{L}_4 q,q\right) \geq C_0 \|q\|^2_{H^1_{per}}.
\end{equation}
Now, we estimate  $\left(\mathcal{L}_3 p,p\right)$. Suppose with out loos of generality that  $\|\phi_\nu\|_{L^2_{per}}=1$. Denote $p_{\perp}=p-p_{\parallel}$, where $p_{\parallel} = (p,\phi_\nu)\phi_\nu$, then from (\ref{equa5}) we obtain that $(p_{\perp},\phi_\nu)=0$ and $(p_{\perp},(\phi_\nu\psi_\nu)')=0$. From  Theorem \ref{teoEstForCua1} it follows that  $(\mathcal{L}_3 p_{\perp},p_{\perp})\geq \widetilde{C_0}\|p_{\perp}\|^2_{L^2_{per}}$.\\

Also consider the normalization $Q_2(u_0)=Q_2(\phi),$ i.e., $\|u_0\|_{L^2_{per}}=\|\phi_\nu\|_{L^2_{per}}$. Then $\|u(t)\|_{L^2_{per}}=1,$ for all  $t\geq 0$, thus $-2(p,\phi_\nu)= \|\xi\|^2_{L^2_{per}}$. Therefore
\[(\mathcal{L}_3p_{\perp},p_{\perp})\geq C_0\|p_{\perp}\|^2_{L^2_{per}} \geq C_0\|p\|^2_{L^2_{per}} -\widetilde{C}_1\|\xi\|^4_{H^1_{per}}.\]
Since  $\left(\mathcal{L}_3\phi_\nu,\phi_\nu\right) <0$, it follows that $\left(\mathcal{L}_3 p_{\parallel},p_{\parallel}\right)\geq -\widetilde{C}_2\|\xi\|^4_{H^1_{per}}$. From the  Cauchy-Schwarz inequality we get that $\left(\mathcal{L}_3 p_{\parallel},p_{\perp}\right)\geq -\widetilde{C}_3\|\xi\|^4_{H^1_{per}}$. Therefore, from the specific form of the operator $\mathcal{L}_3$ we conclude that
\begin{equation}\label{equa9}
\left(\mathcal{L}_3 p,p\right)\geq D_1\|p\|_{H^1_{per}}^2- D_2\|p\|_{H^1_{per}}^3-D_3\|p\|_{H^1_{per}}^4,
\end{equation}
where $D_j >0$ for $j=1,2,3.$\\

Now, using (\ref{equa55}), (\ref{equa8}) and (\ref{equa9}) we arrive at 
\[\Delta\mathcal{B}(t)\geq d_1\|\xi\|_{1,\nu}^2- d_2\|\xi\|_{1,\nu}^3- d_3\|\xi\|_{1,\nu}^4\]
where $d_i>0$, for $i=1,2,3$ and $\|f\|_{1,\nu}^2:=\|f'\|_{L^2_{per}}^2+ \nu\|f\|_{L^2_{per}}^2$.
Using a similar argument as in Benjamin \cite{benjamin1}, we obtain that for any  $\epsilon>0,$ there exists $\delta(\epsilon)>0$ such that if
\[\|u_0-\widetilde{\phi}\|_{1,\nu}<\delta,\ \ \|V_0-\varphi\|_{L^2_{per}}<\delta \ \ \text{and}\ \ \|v_0-\psi\|_{L^2_{per}}<\delta, \]
then
\begin{equation}\label{desiFinalZak}
\rho_{\nu}(u(t),\phi_{\nu})^2=\|\xi(t)\|^2_{1,\nu}<\epsilon,
\end{equation}
for all  $t\geq 0.$ Therefore we obtain the inequality  (\ref{desiImp2}).\\

Finally, using (\ref{equa55}) and the analysis made above for  $\xi$ we obtain that
\[\int_0^L \left[\sqrt{1-c^2}\gamma+\frac{2\phi p}{\sqrt{1-c^2}}+\frac{p^2+q^2}{\sqrt{1-c^2}}\right]^2dx\leq\epsilon\ \ \ \text{and}\ \ \ \int_0^L(c\gamma-\eta)^2 dx\leq\epsilon.\]
Using the two last inequalities, the Cauchy-Schwarz inequality and  (\ref{desiFinalZak}) we arrive at (\ref{desiImp1}), which proves that $(\psi,\varphi,\widetilde{\phi})$ is stable with respect to small perturbations that preserves the $L^2_{per}([0,L])$ norm of  $\widetilde{\phi}.$  The general case follows from the continuity of the map
\[\nu\in\left(\frac{2\pi^2}{L^2},+\infty\right)\mapsto (\psi, \varphi,\widetilde{\phi}).\]
This finishes the proof of the theorem.
\endproof
\section{Stability for the Solitary Wave Solutions}
In this section we improve the result established by Ya Ping in \cite{Wu1}, to obtain a correct  stability result for the solitary wave solutions associated to the Zakharov system. With regard to the Cauchy problem associate to the system (\ref{equaZakha}) we address the reader to the work of Colliander in \cite{Colliander2}, where is obtained a result of global well-posedness  for initial data $(u,v,v_t)(0)\in H^1(\mathbb{R})\times L^2(\mathbb{R})\times \widehat{H}^{-1}(\mathbb{R})$.\\

If $v_t(x,0)\in \widehat{H}^{-1}(\mathbb{R})$ we can rewrite (\ref{equaZakha}) as the equivalent system 
 \begin{equation}\label{Zaksolit}
\left \{
\begin{aligned}
v_t &=-V_x, \\
 V_t &= -(v+|u|^2)_x \\
iu_t &+ u_{xx} =uv.
\end{aligned} \right.
\end{equation}
The solitary wave solutions for this system are given by 
 \[
\psi_{\omega,c}(\xi)=\left(2\omega+\frac{c^2}{2}\right)\ \text{sech}^2\left(\frac{\sqrt{-4\omega-c^2}}{2}\xi\right), \ \ \ \phi_{\omega,c}(\xi)=\sqrt{\frac{(-4\omega-c^2)(1-c^2)}{2}}\ \text{sech}\left(\frac{\sqrt{-4\omega-c^2}}{2}\xi\right)
\]
\[
\text{and}\ \ \ \varphi_{\omega,c}(\xi)=c\left(2\omega+\frac{c^2}{2}\right)\ \text{sech}^2\left(\frac{\sqrt{-4\omega-c^2}}{2}\xi\right).
\]
As in the periodic case, we have to study the spectral properties of the operators 
\begin{equation}
L_{3}=-\frac{d^2}{{dx}^2}-\left(\omega+\frac{c^2}{4}\right)+3\psi\ \ \ \text{and}\ \ \ \
L_{4}=-\frac{d^2}{{dx}^2}-\left(\omega+\frac{c^2}{4}\right)+\psi.
\end{equation}
The spectral properties necessary to obtain our result of stability were established by Ya Ping in \cite{Wu1}. See also \cite{Wu1} for the definition of orbital stability for the solitary wave solutions associated to the Zakharon system.\\

The next theorem is the principal result of this section.
\begin{theo}
Assume that $4\omega+c^2\leq 0$ and $1-c^2>0.$ Then the solitary wave solutions $(\psi_{\omega,c}, \varphi_{\omega,c}, \widetilde{\phi}_{\omega,c}),$ with $ \widetilde{\phi}_{\omega,c}(x)=e^{i\frac{c}{2}x}\phi_{\omega,c}(x),$ are orbitally stable in $X= L^2_{per}([0,L]) \times L^2_{per}([0,L]) \times H^1_{per}([0,L])$ by the flow generated for the system (\ref{Zaksolit}).
\end{theo}
 \proof
 The proof follows the same ideas of Theorem \ref{teoEstForCua2}. We only observe that in this case
\begin{align*}
\Delta\mathcal{B}(t)& =\left(L_3p,p\right)+\left(L_4q,q\right) + \frac{1}{2}\int_{\mathbb{R}}\left[\sqrt{1-c^2}\gamma+\frac{2\phi p}{\sqrt{1-c^2}}+\frac{p^2+q^2}{\sqrt{1-c^2}}\right]^2 dx \\
&+\frac{1}{2}\int_{\mathbb{R}}(c\gamma-\eta)^2 dx-\int _{\mathbb{R}}\frac{4\phi p(p^2+q^2)}{1-c^2}+\frac{(p^2+q^2)^2}{1-c^2}dx. 
\end{align*}
Note that the constant $d_0$ does not appear because the decaying properties of the solitary wave solutions imply that $d_0=0.$ The rest of the proof follows similarly as the result obtained on the periodic case.
 \endproof


\end{document}